\documentclass[final,1p,times,authoryear]{elsarticle}

\usepackage{amsmath, amssymb, amsthm} %

\usepackage{tikz}
\usetikzlibrary{decorations.markings}
\usetikzlibrary{arrows}
\usetikzlibrary{calc,matrix,fit}
\usetikzlibrary{shapes}
\usepackage[all]{xy}
\usepackage{tikz-cd}
\usepackage{amsfonts}
\usepackage{mathtools}
\usepackage{enumerate}
\usepackage{hyperref}
\usepackage{amsmath}
\usepackage{flushend}

\usepackage{listings}
\usepackage[ruled,noend]{algorithm2e}

\newcommand{\Z}{\mathbb{Z}}
\newcommand{\Zset}{\mathbb{Z}}

\newcommand{\Nset}{\mathbb{N}}
\newcommand{\Rset}{\mathbb{R}}

\DeclareMathOperator{\id}{Id}

\renewcommand{\=}{\coloneqq}

\def\Betti(#1,#2,#3) {\beta_{#1}^{#2 , #3}}

\newcommand{\rrdc}{\mbox{\,\(\Rightarrow\hspace{-9pt}\Rightarrow\)\,}}
\newcommand{\lrdc}{\mbox{\,\(\Leftarrow\hspace{-9pt}\Leftarrow\)\,}}
\newcommand{\lrrdc}{\mbox{\,\(\Leftarrow\hspace{-9pt}\Leftarrow\hspace{-5pt}\Rightarrow\hspace{-9pt}\Rightarrow\)\,}} 

\newtheorem{remark}{Remark}
\newtheorem{example}{Example}

\newcommand\TODO[3]{\hbox to 0pt{\textcolor{#1}{$^\bullet$}}\marginpar{\textcolor{#1}{#2: #3}}}

\definecolor{dblackcolor}{rgb}{0.0,0.0,0.0}
\definecolor{dbluecolor}{rgb}{.01,.02,0.7}
\definecolor{dredcolor}{rgb}{0.8,0,0}
\definecolor{dgraycolor}{rgb}{0.30,0.3,0.30}

\usepackage{mparhack}

\newcommand{\sagemath}{\text{SageMath}}
\newcommand{\kenzo}{\text{Kenzo}}

\newtheorem{theorem}{Theorem}[section]
\newtheorem{lemma}[theorem]{Lemma}
\newtheorem{proposition}[theorem]{Proposition}
\newtheorem{corollary}[theorem]{Corollary}

\theoremstyle{definition}
\newtheorem{definition}[theorem]{Definition}

\begin{document}

\title{Computing the homology of universal covers via effective homology and discrete vector fields \footnote{The first author has been partially supported by PID2020-114750GB-C31, E22\_20R: Álgebra y Geometría, and by project Inicia 2023/02 funded by La Rioja Government (Spain). The second 
author is supported by grant PID2020-116641GB-I00 funded by MICIU/ AEI/ 10.13039/501100011033.}}

\author[1]{Miguel A. Marco-Buzun\'ariz}
\ead{mmarco@unizar.es}

\author[2]{Ana Romero}
\ead{ana.romero@unirioja.es}


\address[1]{University of Zaragoza, Spain }

\address[2]{University of La Rioja, Spain }


\begin{abstract}
Effective homology techniques allow us to compute homology groups of a wide family of topological spaces. By the Whitehead tower method, this can also be used to compute higher homotopy groups. However, some of these techniques (in particular, the Whitehead tower) rely on the assumption that the starting space is simply connected. For some applications, this problem could be circumvented by replacing the space by its
universal cover, which is a simply connected space that
shares the higher homotopy groups of the initial space.
In this paper, we formalize a simplicial construction for the universal cover, and represent it as a twisted
cartesian product.

As we show with some examples, the universal cover of a space with effective homology does not necessarily have effective homology in general. We show two independent sufficient conditions that can ensure it: one is based on a nilpotency property of the fundamental group, and the other one on discrete vector fields.

Some examples showing our implementation of these constructions in both \sagemath\ and \kenzo\ are shown, together with an approach to compute the homology of the universal cover when the group is abelian even in some cases where there is no effective homology, using the twisted homology of the space.
\end{abstract}


\maketitle

\section{Introduction}


The \emph{effective homology} method, implemented in the computer algebra system~\citep{Kenzo},  was developed by Francis Sergeraert in~\citep{Ser94} with the aim of calculating homology and homotopy groups of complicated topological spaces, allowing in particular to work with spaces of infinite type. Topological spaces are represented in Kenzo by means of \emph{simplicial sets}, a combinatorial structure that generalizes the notion of simplicial complex. The computation of homology groups is done in Kenzo by means of chain equivalences between the initial simplicial set and a finite type chain complex in which the homology can be computed by matrix diagonalization operations. For the computation of homotopy groups, the Whitehead tower method~\citep{Whi52} is used.

Most algorithms in Kenzo and, in general, in the effective homology method, require the simplicial sets to be $1$-reduced (that is, with only one vertex and without non-degenerate edges). In particular, this condition is necessary for the application of the Kenzo's implementation of the Whitehead tower method. In a previous work~\citep{mathematics-2021}, a Kenzo interface and an optional package for SageMath~\citep{sagemath} were developed. In that work, we also developed and integrated an algorithm to compute homotopy groups of simply connected simplicial sets that are not necessarily $1$-reduced. However, the condition that the simplicial set be simply connected  is a necessary condition for the Whitehead tower method~\citep{Whi52}.

A \emph{cover} of a connected topological space $X$ is a topological space $Y$ with a map $f:Y\rightarrow X$ such that for every $x \in X$, there exists an open neighborhood $U$ of $x$ such that $f^{-1}(U)$ is a disjoint union of homeomorphic copies of $U$. If $Y$ is simply connected, then $Y$ is said to be a \emph{universal cover} of $Y$, and satisfies $\pi_i(Y) \cong \pi_i(X)$ for all $i \geq 2$. The universal cover of a topological space $X$ can be considered as an approach to study $X$ by lifting paths and other geometric objects to the simpler universal cover. This simplifies the resolution of certain topological and geometric issues on $X$, as they can be streamlined to the examination of the universal cover.  In particular, the universal cover can be useful in applying the Whitehead tower method.

In this work, we present algorithms for computing a simplicial model of the universal cover of a space, represented as a twisted cartesian product, and its effective homology by using homological perturbation theory. The effective homology of the universal cover can be determined with our algorithms when the effective homology of the initial space satisfies one of two independent conditions: one is based on a nilpotency property of the fundamental group, and the other one on discrete vector fields. We also present an approach to compute the homology of the universal cover when the fundamental group is abelian even in some cases where there is no effective homology, using the twisted homology of the space.
 Our algorithms have been implemented in the computer algebra systems \sagemath~\citep{sagemath} and \kenzo\, using the Kenzo interface and the optional package for SageMath that we developed in~\citep{mathematics-2021}.

This work presents a revised and extended version of our previous conference paper~\citep{MRD23}. In addition to presenting more details, with the necessary preliminaries and proofs of all of our results in Sections~\ref{sec:simplicial_universal_cover} and~\ref{sec:effective_homology_algorithm} (that in~\citep{MRD23} were only briefly sketched), we have included two new sections with new results. On the one hand, we have studied now the behavior of some algebraic topology constructors and tools with respect to our construction of the universal cover and its effective homology, considering in particular cartesian products and discrete vector fields. On the other hand, we have also examined the case of abelian infinite fundamental group, when the universal cover might not have effective homology.

The paper is structured as follows. In the next section, we present preliminary definitions and results about simplicial sets, effective homology and universal covers. In Section~\ref{sec:simplicial_universal_cover}, our simplicial construction of the universal cover as a twisted cartesian product is described. Section~\ref{sec:effective_homology_algorithm} shows the  algorithms for computing the effective homology of universal covers, making use of homological perturbation theory. Then, we explain the implementation of the algorithms and show some didactic examples in Section~\ref{sec:implementation_and_examples}. Cartesian products and discrete vector fields and their relation with universal covers and their effective homology is studied in Section~\ref{sec:constructors_and_tools}, and the case of abelian infinite fundamental group is considered in Section~\ref{sec:infinite_fundamental_group}. Finally, we describe the conclusions and possible further work in Section~\ref{sec:conclusions}.

\section{Preliminaries}
\label{sec:prelim}

\subsection{Simplicial sets}
\label{subsec:simplicial_top}

In this section, we introduce the definition and some basic constructions of simplicial sets, following \citep{May67}.

\begin{definition}
\label{defn:smob}
Let $\mathcal{D}$ be a category. The category $s\mathcal{D}$ of \emph{simplicial objects in} $\mathcal{D}$ is defined as follows.
An object $X \in s\mathcal{D}$ consists of 
\begin{itemize}
\item for each integer $n \geq 0$, an object $X_n \in \mathcal{D}$;
\item for every pair of integers $(i,n)$ such that $0 \leq i \leq n$, \emph{face} and \emph{degeneracy} maps $\partial_i: X_n \rightarrow X_{n-1}$ and $s_i: X_n \rightarrow X_{n+1}$ (which are morphisms in the category~$\mathcal{D}$) satisfying the \emph{simplicial identities}:

\begin{align*}
\partial_i \partial_j &= \partial_{j-1} \partial_i && \text{if } i < j \\
s_i s_j &= s_{j+1} s_i && \text{if } i \leq j \\
  & && \text{if } i < j \\
  \partial_i s_j &=
    \smash{\left\{\begin{array}{@{}l@{}}
      s_{j-1} \partial_i \\[\jot]
      \mathrm{Id} \\[\jot]
      s_j \partial_{i-1}
    \end{array}\right.} && \text{if } i=j, j+1 \\
  & && \text{if } i > j+1
\end{align*}
\end{itemize}
Let $X$ and $Y$ be simplicial objects. A \emph{simplicial map} (or \emph{simplicial morphism}) $f: X \rightarrow Y$ consists of maps $f_n: X_n \rightarrow Y_n$ (which are morphisms in $\mathcal{D}$) that commute with the face and degeneracy operators, that is
$f_{n-1} \partial_i = \partial_i  f_n$ and $f_{n+1} s_i = s_i f_n$ for all $0 \leq i \leq n$.
\end{definition}

If $\mathcal{D}$ is a subcategory of sets (which we will assume from now on), the elements of $X_n$ are called the \emph{$n$-simplices} of $X$.

\begin{definition}
\label{def:deg_nondeg}
An $n$-simplex $x\in X_n$ is \emph{degenerate} if $x=s_jy$ for some $y \in X_{n-1}$ and some $0 \leq j < n$; otherwise, $x$ is called \emph{non-degenerate}. 
We denote by $X^{D}_n$ the set of degenerate $n$-simplices and by $X^{ND}_n$ the set of non-degenerate $n$-simplices of $X$.
\end{definition}

A \emph{simplicial set} is a simplicial object in the category of sets.
A \emph{simplicial group} $G$ is a simplicial object in the category of groups; in other words, it is a simplicial set where each $G_n$ is a group and the face and degeneracy operators are group morphisms.

A simplicial set $X$ has a canonically associated chain complex $C_*(X)=(C_n(X),d_n)$, where each chain group $C_n(X)$ is defined as the free $\Z$-module generated by $X_n$, and the differential $d_n : C_n(X) \to C_{n-1}(X)$ is defined as the alternating sum of faces, $d_n \= \sum_{i=0}^n (-1)^i \partial_i$. In this work we will assume all the chain complexes associated with simplicial sets to be \emph{normalized} (see \cite[Ch.~5]{May67}), which intuitively means that only non-degenerate simplices are considered as generators of the chain groups.

\begin{definition}
\label{def:Cartesian_KL}
The \emph{cartesian product} $X\times Y$ of two simplicial sets $X$ and $Y$ is the simplicial set whose set of $n$-simplices is $(X\times Y)_n \= X_n \times Y_n$, with coordinate-wise defined face and degeneracy maps: if $(x,y)\in (X\times Y)_n$, then
\begin{align*}
\partial_i (x,y) &\= (\partial_i x , \partial_i y), \qquad 0\le i\le n ; \\
s_i (x,y) &\= (s_i x , s_i y), \qquad 0\le i\le n .
\end{align*}
\end{definition}


\begin{definition}
\label{def:Twisted_GB}
A \emph{twisting operator} from a simplicial set $B$ to a simplicial group $G$ is a map $\tau : B\to G$ of degree $-1$, that is a collection of maps $\tau = \{ \tau_n : B_n \to G_{n-1} \}_{n\ge 1}$, satisfying the following identities, for any $n\ge 1$ and for any $b\in B_n$:
\begin{align*}
\partial_i ( \tau b) &= \tau (\partial_i b), \qquad 0\le i<n-1 , \\
\partial_{n-1} ( \tau b) &= \tau (\partial_n b)^{-1} \cdot \tau (\partial_{n-1} b), \\
s_i (\tau b) &= \tau (s_i b), \qquad 0\le i\le n-1 , \\
e_n &= \tau (s_n b) ,
\end{align*}
where $e_n$ is the identity element of $G_n$.
\end{definition}

We defined twisting operators in a slightly different (but equivalent) way from \citep{May67}, in order to agree with the definition implemented in the Kenzo system.

\begin{definition}
\label{def:twcrpr}
Given a simplicial group $G$, a simplicial set $B$ and a twisting operator $\tau : B\to G$,
the \emph{twisted (cartesian) product} $E(\tau)\= G \times_{\tau} B$ is the simplicial set whose set of $n$-simplices is $E(\tau )_n = (G\times_{\tau} B)_n \= G_n \times B_n$ and whose face and degeneracy maps are defined in the following way: if $(g,b)\in (G\times_{\tau} B)_n$, then
\begin{align*}
\partial_i (g,b) &\= (\partial_i g , \partial_i b), \qquad 0\le i < n , \\
\partial_n (g,b) &\= (\tau (b) \cdot \partial_n g , \partial_n b) ; \\
s_i (g,b) &\= (s_i g , s_i b), \qquad 0\le i\le n .
\end{align*}
\end{definition}

It can be easily shown that the identities defining a twisting operator $\tau$ are equivalent to the simplicial identities of $G\times_{\tau} B$. 


\subsection{Effective homology and homological perturbation theory}
\label{subsec:effective_homology}

We present now the main definitions and ideas of the effective homology method, introduced in \citep{Ser94} and explained in depth in~\citep{RS02} and~\citep{RS06}. We will also recall some results in homological perturbation theory that we will need in our work.

\begin{definition}
\label{def:red}
A \emph{reduction} $\rho\equiv(D_\ast \rrdc C_\ast)$ between two
chain complexes $D_\ast$ and $C_\ast$ is a triple $(f,g,h)$ where:
\begin{enumerate}
        \item  The components
$f$ and $g$ are chain complex morphisms $f: D_\ast \rightarrow C_\ast$ and $g: C_\ast \rightarrow D_\ast$;
\item
The component $h$ is a homotopy operator $h:D_\ast\rightarrow D_{\ast+1}$ (a graded group homomorphism of degree +1);
\item The following relations are satisfied:
  \begin{enumerate}
      \item    $f  g = \mbox{id}_{C_\ast}$;
      \item   $g f + d_{D_\ast} h + h  d_{D_\ast}
        = \mbox{id}_{D_\ast}$;
  \item  {$f  h = 0$;} $h   g = 0$; $h   h = 0$.
  \end{enumerate} \end{enumerate}
\end{definition}

\begin{remark}
These relations show that $D_\ast$ is the direct sum of $C_\ast$ and a contractible (acyclic) complex. 
In particular, this implies that the graded homology groups \(H_\ast(D_\ast)\) and \(H_\ast(C_\ast)\) to be canonically isomorphic.
\end{remark}

\begin{definition}
A \emph{(strong chain) equivalence} $\varepsilon \equiv (C_\ast \lrrdc E_\ast)$ between two complexes $C_\ast$ and~$E_\ast$ is a triple $(D_\ast,\rho,\rho ')$ where $D_\ast$ is a chain complex and $\rho$ and $\rho'$ are reductions from $D_\ast$ over $C_\ast$ and  $E_\ast$ respectively: $C_\ast \stackrel{\rho}{\lrdc} D_\ast \stackrel{\rho'}{\rrdc} E_\ast.$
\end{definition}

\begin{definition}
An \emph{effective} chain complex $C_\ast$ is a free chain complex (i.e., a chain complex consisting of free $\Zset$-modules) where each group $C_n$ is finitely generated, and there is an algorithm that returns a $\Zset$-base $\beta_n$ in each degree $n$ (for details, see \cite{RS02}).
\end{definition}

The homology groups of an effective chain complex $C_\ast$ can easily be determined by means of some diagonalization algorithms on matrices (see \cite{KMM04}).

\begin{definition}
An \emph{object with effective homology} is a triple $(X,EC_\ast,\varepsilon)$ where $X$ is an object (e.g. a simplicial set, a topological space) possessing a canonically associated free chain complex $C_\ast(X)$, $EC_\ast$ is an effective chain complex and $\varepsilon$ is an equivalence between $C_\ast(X)$ and $EC_\ast$, $C_\ast(X) \stackrel{\varepsilon}{\lrrdc} EC_\ast$.
\end{definition}


The notion of object with effective homology makes it possible in this way to
compute homology groups of complicated spaces by using effective complexes and computing their homology groups. The effective homology method is based on the following idea: Given some topological
spaces $X_1, \ldots, X_n$, a topological constructor $\Phi$ produces a new topological space
$X$. If effective homology versions of the spaces $X_1, \ldots , X_n$ are known, then an
effective homology version of the space $X$ can also be built, and this version
allows us to compute the homology groups of $X$, even if it is not of finite type. In this way, we can obtain the effective homology of objects that are constructed using these constructors and other objects with effective homology. This method has been implemented in the Kenzo system \citep{Kenzo}, a Common Lisp program devoted to Symbolic Computation in Algebraic Topology, which has made it possible to determine homology and homotopy groups of complicated spaces. 
The computation of homotopy groups of simply connected simplicial sets with effective homology is dealt with in Kenzo by means of the Whitehead tower method~\citep{Whi52} and requires the construction of a sequence of fibrations involving spaces of infinite types (following the algorithm in \citep{Rea96}).

Now, we introduce the main results in homological perturbation theory, which describe how a perturbation (a modification of the differential of a chain complex) transmits through a reduction.

\begin{definition}
Let $C_*=(C_n,d_n)_{n \in \Zset}$ be a chain complex. A \emph{perturbation} $\delta$ of the differential $d$ is a family of morphisms $\delta=\{ \delta_n: C_n \to C_{n-1}\}_{n \in \Zset}$ such that the sum $d + \delta$ is again a differential, that is $(d + \delta)^2=0$ holds (meaning $(d_{n-1} + \delta_{n-1}) (d_n +\delta_n)= 0 $, for all $ n\in \Zset$).
\end{definition}

We call $C'_*=(C_n,d_n+\delta_n)_{n \in \Zset}$ the \emph{perturbed} chain complex obtained from~$C_*$ by introducing the perturbation $\delta$.

\begin{theorem} [Trivial Perturbation Lemma, TPL] 
\label{thm:tpl}
Let $C_*= (C_n, d_{C_n})_{n \in \Zset}$ and $D_*=(D_n,d_{D_n})_{n \in \Zset}$ be two chain complexes, $\rho=(f,g,h): C_* \rrdc D_*$ a reduction, and $\delta_D$ a perturbation of the differential
$d_D$. Then a reduction $\rho'=(f',g',h'): C'_* \rrdc D'_*$ can be built, where:
\begin{enumerate}
\item$C'_* = (C_* , d_{C}+ g  \delta_D  f)$ is the perturbed chain complex obtained from $C_*$ by introducing the perturbation $g  \delta_D  f$;
\item $D'_* = (D_* , d_D +\delta_D )$ is the perturbed chain complex obtained from $D_*$ by introducing the perturbation $  \delta_D$;
\item the maps of the new reduction $\rho'=(f',g',h')$ are given by $f'\= f$, $g'\= g$, $h'\= h$.
\end{enumerate}
\end{theorem}

\begin{theorem}[Basic Perturbation Lemma, BPL,~\cite{Brown1967}]
\label{thm:bpl}
Let $C_*=(C_n,d_{C_n})_{n \in \Zset}$ and $D_*=(D_n,d_{D_n})_{n \in \Zset}$ be two chain complexes, $\rho=(f,g,h): C_* \rrdc D_*$ a reduction, and $\delta_C$ a perturbation of the differential
$d_C$. Suppose that the composition $h  \delta_C$ satisfies the following \emph{nilpotency condition}: for every $x \in C_*$ there exists a non-negative integer $m=m(x) \in \Nset$ such that
$(h \delta_C)^m(x)=0$. Then a reduction
$\rho'=(f',g',h'): C'_* \rrdc D'_*$ can be built, where:
\begin{enumerate}
\item $C'_*= (C_* , d_{C}+  \delta_C )$ is the perturbed chain complex obtained from $C_*$ by introducing the perturbation $  \delta_C$;
\item $D'_*= (D_* , d_{D}+  \delta_D )$ is the perturbed chain complex obtained from $D_*$ by introducing the perturbation $\delta_D \= f  {\delta_C}  \varphi  g =
               f  \psi  {\delta_C}  g$;
\item the maps of the new reduction $\rho'=(f',g',h')$ are given by
\begin{equation*}
f'  \= f  \psi , \qquad
g'  \= \varphi  g , \qquad
h'  \= \varphi  h = h  \psi ,
\end{equation*}
\end{enumerate}
with the operators $\varphi$ and $\psi$ given by
\begin{equation*}
\varphi  \=\sum_{i=0}^{\infty}{(-1)^i(h \delta_C)^i} , \qquad
\psi  \=\sum_{i=0}^{\infty}{(-1)^i(\delta_C  h)^i} ,
\end{equation*}
the convergence of these series being guaranteed by the nilpotency condition.
\end{theorem}

Given a chain complex $C_*$, the \emph{trivial reduction} $\id = (f,g,h) :C_* \rrdc C_*$ is the reduction with $f=g=\id$ and $h=0$. To finish this subsection, we state a simple result (see \cite[Ch.~5]{RS06}) that describes the behavior of reductions with respect to the tensor product.

\begin{proposition}
\label{prop:red_tensor}
Let $\rho=(f,g,h): C_* \rrdc D_*$ and $\rho'=(f',g',h'): C'_* \rrdc D'_*$ be two reductions. Then a reduction $\rho''\=\rho \otimes \rho' =(f'',g'',h''): C_* \otimes C'_* \rrdc D_* \otimes D'_*$ is given by:
\[ f'' \= f \otimes f' , \qquad g'' \= g \otimes g'  , \qquad h'' \= h \otimes \id_{C'_\ast} + (g  f) \otimes h'.
\]
\end{proposition}

\subsection{Universal covers}
\label{sec:universal_covers}

For completeness, we include here a brief summary of universal covers and their properties. See, for example, \citep{hatcher2002algebraic} for a fully detailed explanation.

A continuous map $f: \tilde X \to X$ between two path-connected spaces is said to be a \emph{covering map} if  for every $x\in X$, there exists a neighborhood $U$ such that $f^{-1}(U)$ is a disjoint union of homeomorphic copies of $U$ (and the restriction of $f$ to each such copy is in fact a homeomorphism). In this situation, $\tilde X$ is said to be a \emph{cover} of $X$.

A relevant property of covering maps is that they have \emph{unique elevation of paths}: given a point $\tilde x_0\in \tilde X$ and a path $\gamma:I \to X$ starting at $x_0=f(\tilde x_0)$, there exists a unique lift $\gamma_{x_0}:I\to \tilde X$ starting at $x_0$. This implies the induced map $f^* : \pi_1(\tilde X)  \to \pi_1(X)$ to be injective, so we can see the fundamental group of $\tilde X$ as a subgroup of the fundamental group of $X$. If this subgroup is normal, the cover is said to be \emph{regular} (or \emph{Galois}). In this case, the base space $X$ is a quotient of $\tilde X$ by the group of so-called \emph{deck transformations}. This group is isomorphic to the quotient of $\pi_1(X)$ by $\pi_1(\tilde X)$.

If $\tilde X$ is simply connected, it is said to be a \emph{universal cover} of $X$. Every connected, locally path connected and semilocally $1$-connected space has a universal cover (and, in fact, only one, in the sense that two such covers would be homeomorphic). That is, such spaces are the quotient of a certain simply connected space by the action of some group.

Since covering maps are locally trivial fibrations, with discrete fiber, we can apply the long exact sequence to get isomorphisms $\pi_i(\tilde X) \to \pi_i(X)$ for $i\geq 2$. That is, the higher homotopy groups of a space coincide with the ones of its universal cover.

As we have said before, one of the possible applications of effective homology is to compute the homotopy groups of spaces through the Whitehead tower method. This requires the spaces to be simply connected, but we can satisfy that condition if we compute the universal cover of the space. In fact, this can be seen as a first step in the Whitehead tower.

\section{A simplicial construction for universal covers}
\label{sec:simplicial_universal_cover}




Given a connected simplicial set $X$, a presentation of its fundamental group can be found as follows:
\begin{enumerate}
    \item Choose a maximal tree $T$ in the $1$-skeleton of $X$.
 \item Take a generator $g_e$ for each edge $e$ that is not in $T$.
 \item  For each non-degenerate $2$-simplex $\sigma$ in $X$, add a relation given by $g_{\partial_2\sigma}g_{\partial_0\sigma}g_{\partial_1\sigma}^{-1}$ (assuming that the edges in $T$ correspond to the trivial element). This relation represents the fact that a closed path that follows the boundary of a triangle can be retracted to a constant path.
\end{enumerate}

 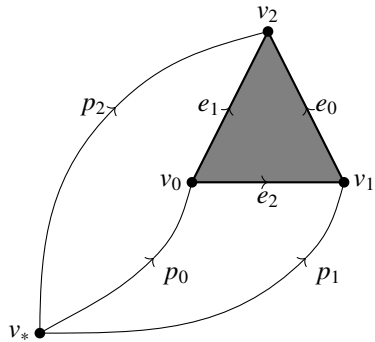
\begin{figure}[h]
 \begin{center}
 \begin{tikzpicture}
     \begin{scope}[scale=2]
     \node[anchor=east] (v0) at (0,0)  {$v_0$};
     \node[anchor=west] (v1) at (1,0)  {$v_1$};
     \node[anchor=south] (v2) at (0.5,1)  {$v_2$};
     \node[anchor=east] (v) at (-1,-1) {$v_*$};
     \node[anchor=north] (e2) at (0.5,0) {$e_2$};
     \node[anchor=west] (e0) at (0.75,0.5) {$e_0$};
     \node[anchor=east] (e1) at (0.25,0.5) {$e_1$};

     \node[anchor=north west] (p0) at (-0.25,-0.5) {$p_0$};
    \node[anchor=north west] (p1) at (0.75,-0.5) {$p_1$};
    \node[anchor=east] (p2) at (-0.5,0.5) {$p_2$};

    \draw[fill=gray,thick] (0,0) -- (1,0) -- (0.5,1) -- cycle;

    \fill (-1,-1) circle[radius=1pt];
    \fill (0,0) circle[radius=1pt];
    \fill (1,0) circle[radius=1pt];
    \fill (0.5,1) circle[radius=1pt];

    \draw[->] (0,0) -- (0.5,0);
    \draw[->] (1,0) -- (0.75,0.5);
    \draw[->] (0,0) -- (0.25,0.5);

    \draw[->,out=30,in=225] (-1,-1) to (-0.25,-0.5);
    \draw[out=45,in=255] (-0.25,-0.5) to (0,0);
    \draw[out=0,in=225,->] (v) to (0.75,-0.5);
    \draw[out=45,in=255] (0.75,-0.5) to (1,0);
    \draw[->,out=90,in=225] (-1,-1) to (-0.5,0.5);
    \draw[out=45,in=195] (-0.5,0.5) to (0.5,1);

 \end{scope}
 \end{tikzpicture}
 \end{center}
 \caption{Relation given by a $2$-simplex}
\end{figure}



The result is a presentation of $G:=\pi_1(X)$. Note that we also get a map $\tau '$ that assigns an element of the group to each edge (assuming that the edges in $T$ correspond to the trivial element). This assignment can be lifted to higher-dimensional simplices by taking the first face. In order to work effectively with this presentation, we need to be able to solve the word problem in this presentation (which cannot be ensured in the general case, but is solved for several families of groups, including the cases of finite, free, abelian, polycyclic or simple groups).

Now consider a multiplicative group $H$ and a surjective group morphism $\chi: \pi_1(|X|,v_*)\to H$.
We define
a map
$$
\tau: X \to H
$$

\noindent given by the following recursive rules:

\begin{itemize}
 \item $\forall v \in X_0,\tau(v)=1$
\item $\forall e \in T, \tau(e)=1$
\item $\forall v \in X_0, \tau(s_0v)= 1$
\item $\forall e\in X_1 \setminus T, \tau(e)=\chi(\tau'(e))$
\item $\forall \sigma \in X_n$  with $n>1$, $\tau(\sigma)=\tau(\partial_0\sigma)$
\end{itemize}


\begin{lemma}
 With the previous assignation rules, for each $\sigma\in X_n$ with $n\geq 2$,
 the following equality holds:
 $$
 \tau(\partial_{n-1}\sigma)= \tau(\partial_n\sigma) \cdot \tau(\partial_0\sigma).
 $$
\end{lemma}

\begin{proof}
 We prove it by induction on $n$.

For $n=2$, it is a direct consequence of the relations in the presentation of the group.

Now assume that it is true for simplices of dimension up to $n-1$, and take
a simplex $\sigma$ of dimension $n$. Applying the recursive definition of
$\tau$, we get:

$$
\begin{array}{rcl}
    \tau(\partial_{n-1}\sigma) & = &  \tau(\partial_0\partial_{n-1}\sigma) = \\
    & = & \tau(\partial_{n-2}\partial_{0}\sigma)= \\
    & = & \tau(\partial_{n-1}\partial_{0}\sigma) \cdot \tau(\partial_0\partial_{0}\sigma)= \\
& = &\tau(\partial_0\partial_{n}\sigma) \cdot \tau(\partial_0\partial_{0}\sigma) = \\
& = & \tau(\partial_n\sigma) \cdot\tau(\partial_0\sigma)
\end{array}
$$

\end{proof}

\begin{lemma}
 For each $\sigma\in X_n$ with $n\geq 2$, and $0\leq i < n-1$,
 $\tau(\sigma)=\tau(\partial_i\sigma)$.
\end{lemma}

\begin{proof}
By definition $\tau(\partial_i\sigma) =\tau(\partial_0\partial_i\sigma)$.
Iterating this $i$ times, we get $\tau(\partial_0\dots \partial_0\partial_i\sigma)$
which, by the properties of the face maps in a simplicial set, equals $\tau(\partial_0\partial_0\dots \partial_0\sigma)$. Again, iterating the
definition $i+1$ times, this is equal to $\tau(\sigma)$.
\end{proof}

We will now construct a new simplicial set $\tilde{X}$ with the following rules:

\begin{itemize}
    \item The sets of simplices are $\tilde{X}_n= H \times X_n$
    \item The degeneracy maps are $\tilde{s}_i(h,\sigma) = (h,s_i\sigma)$
    \item The face maps for an $n$-dimensional simplex $(h,\sigma)$ are given by:
    \begin{itemize}
        \item $\tilde{\partial}_i(h,\sigma) = (h,\partial_i\sigma)$, for $i<n$
        \item $\tilde{\partial}_n(h,\sigma) = (h \cdot \tau(\sigma)^{-1}, \partial_n\sigma)$
    \end{itemize}

\end{itemize}

\begin{lemma}
    The maps $\{\tilde{s}_i\}$ and $\{\tilde{\partial}_i\}$ define a
    structure of simplicial set in $\tilde{X}$.
\end{lemma}

\begin{proof}
The properties of the compositions of the degeneracies are a consequence of the ones in $X$.
To check the rule for composition of faces, we distinguish cases for an $n$-dimensional
simplex. 

For $n\geq 2$, we have

$$
\begin{array}{c}
\tilde{\partial}_{n-1}\tilde{\partial}_n(h,\sigma)=\tilde{\partial}_{n-1}(h \cdot \tau(\sigma)^{-1},\partial_n\sigma )= \\
    = (h \cdot \tau(\sigma)^{-1} \cdot \tau(\partial_{n}\sigma)^{-1},\partial_{n-1}\partial_n\sigma ) = \\
    = (h \cdot  \tau(\partial_0\sigma)^{-1} \cdot \tau(\partial_n\sigma)^{-1}, \partial_{n-1}\partial_{n-1}\sigma)= \\
    = (h \cdot \tau(\partial_{n-1}\sigma)^{-1},\partial_{n-1}\partial_{n-1}\sigma) = \\
    \tilde{\partial}_{n-1}(h,\partial_{n-1}\sigma)=\tilde{\partial}_{n-1}\tilde{\partial}_{n-1}(h,\sigma)
\end{array}.
$$

If $i<n-1$,

$$
\begin{array}{c}
\tilde{\partial}_{i}\tilde{\partial}_n(h,\sigma)=\tilde{\partial}_{i}(h \cdot \tau(\sigma)^{-1},\partial_n\sigma)= \\
    = (h \cdot \tau(\sigma)^{-1},\partial_{i}\partial_n\sigma) = \\
    = (h \cdot \tau(\sigma)^{-1},\partial_{n-1}\partial_i\sigma) = \\
    = (h \cdot \tau(\partial_i\sigma)^{-1},\partial_{n-1}\partial_i\sigma) = \\
\tilde{\partial}_{n-1}\tilde{\partial}_i(h,\sigma)
\end{array}.
$$

For $i<j<n$, the relations between $\partial_i$ and $\partial_j$ transfer directly
to $\tilde\partial_i$ and $\tilde\partial_j$.
\end{proof}

As a direct consequence of the definition of $\tilde\partial_i$, we get the following result.

\begin{lemma}
 The map
 $$
 \begin{array}{rcl}
 \pi: \tilde{X} & \to &  X \\
 (h,\sigma) & \mapsto & \sigma
\end{array}
 $$

 is a simplicial map.

\end{lemma}

This implies that $\pi$ induces a continuous map $|\tilde{X}| \to   |X|$.

We will now see that it is indeed a covering map. To do so, we need some preliminary results.
\begin{lemma}
    \label{lemma:unica-elevacion-incidencia}

Let $\sigma_0,\sigma_1$ be simplices in $X$ such that $\partial_{a_1}\cdots \partial_{a_m}\sigma_1=\sigma_0$ for some $a_1<\cdots<a_m$. Then for every $g\in H$, there exists a unique $g'\in H$ such that $\tilde\partial_{a_1} \cdots \tilde\partial_{a_m}(g',\sigma_1)=(g,\sigma_0)$.
\end{lemma}
\begin{proof}
Applying the rules that define $\tilde\partial_i$, we get that
$\tilde\partial_{a_1} \cdots \tilde\partial_{a_m}(1,\sigma_1)=
(h, \partial_{a_1}\cdots \partial_{a_m}\sigma_1)=(h,\sigma_0) $ for some $h\in H$. Then, the equation
$$
\tilde\partial_{a_1} \cdots \tilde\partial_{a_m}(g',\sigma_1)=(g,\sigma_0)
$$
is equivalent to
$$
(g'h,\sigma_0)=(g,\sigma_0)
$$
and the only solution is $g'=gh^{-1}$.
\end{proof}

Notice that, by the way the degeneration maps are defined in $\tilde X$, this result can be trivially generalized to the case where we have a sequence of both faces and degenerations.

In the situation of the previous lemma, we will say that $\sigma_1$ is \emph{incident} to $\sigma_0$.





\begin{lemma}
\label{lemma:caras-degeneradas-pegan}
Let $\sigma,\sigma_0$ be non-degenerate incident simplices such that $\partial_{a_1}\cdots \partial_{a_m}\sigma=s_{b_1}\cdots s_{b_{l+1}}\sigma_0$, and $\partial_i\partial_{a_1}\cdots \partial_{a_m}\sigma=s_{c_1}\cdots s_{c_l}\sigma_0$ for some $i$, $a_1\ldots,a_m$, $b_1\ldots,b_{l+1}$ and some $c_1,\ldots,c_{l}$ (that is, there is a face $F$ of $\sigma$ that is a degeneration of $\sigma_0$, and so is a face of $F$).

If $\tilde\partial_{a_1}\cdots \tilde\partial_{a_m}(g,\sigma)= \tilde s_{b_1}\cdots \tilde s_{b_{l+1}}(h,\sigma_0)$ and $\tilde\partial_i\tilde\partial_{a_1}\cdots \tilde\partial_{a_m}(g',\sigma)=\tilde s_{c_1}\cdots\tilde s_{c_l}(h,\sigma_0)$, then $g$ and $g'$ must be equal.
\end{lemma}

\begin{proof}
Let be $d$ the dimension of $\sigma_0$. If $i<d+l+1$, then clearly

\begin{align*}
    \tilde\partial_i\tilde\partial_{a_1}\cdots \tilde\partial_{a_m}(g',\sigma)  &= \tilde \partial_i (g'h',\partial _{a_1}\cdots \partial_{a_m}\sigma) \\
    &= (g'h', \partial_i\partial _{a_1}\cdots \partial_{a_m}\sigma)
\end{align*}
for some $h'$ that depends only on $a_1\ldots,a_m$, $\sigma$ and $\tau$.

whereas

$$
\tilde\partial_{a_1}\cdots \tilde\partial_{a_m}(g,\sigma) =
(gh', \partial _{a_1}\cdots \partial_{a_m}\sigma).
$$

with the same $h'$ as before.
If both are degenerations of the same simplex, we have that $g'=g$ (since the degeneration maps in $\tilde X$ don't affect the first coordinate).

Let us consider the case where $i = d + l +1 $. Without loss of generality, we can assume that $b_1 > \cdots > b_{l+1}$. Looking at the dimensions of the simplices, necessarily $b_i \leq  d + l +1 - i$; so $b_1 \leq d + l  $.

If $b_1 < d+l$, we have that

\begin{align*}
    \partial_{d+l+1}s_{b_1}\cdots s_{b_{l+1}} \sigma_0 &= s_{b_1}\partial_{d+l}s_{b_2}\cdots s_{b_{l+1}}\sigma_0 \\
    &= \cdots \\
    &= s_{b_1}\cdots s_{b_{l+1}}\partial_{d}\sigma_0
\end{align*}
but that would imply that $\partial_i\partial_{a_1}\cdots \partial_{a_m}\sigma$ is both a degeneration of $\sigma_0$ (that is a non-degenerate simplex of dimension $d$), and of $\partial_{d}\sigma_0$ (that is a $d-1$ dimensional simplex). Since this cannot happen, necessarily $b_1$ must be equal to $d+l$.

In that case,

\begin{align*}
    \tilde\partial_{d+l+1}\tilde\partial_{a_1}\cdots \tilde\partial_{a_m}(g',\sigma) &= \tilde\partial_{d+l+1} \tilde s_{d+l}s_{b_2}\cdots \tilde s_{b_{l+1}}(h,\sigma_0) \\
    &= s_{b_2}\cdots \tilde s_{b_{l+1}}(h,\sigma_0) \\
    &= \tilde\partial_{d+l} \tilde s_{d+l}s_{b_2}\cdots \tilde s_{b_{l+1}}(h,\sigma_0) \\
    &= \tilde\partial_{d+l}\tilde\partial_{a_1}\cdots \tilde\partial_{a_m}(g',\sigma)
\end{align*}

and we can apply the previous case.

\end{proof}

\begin{proposition}
  \label{lemma:covering_map}
 The map induced by $\pi$ is a covering map.
\end{proposition}

\begin{proof}

Let $x \in |X|$, consider the smallest $n$ such that $x$ lives in the
$n$-skeleton of $|X|$. There is a unique non-degenerate $n$-dimensional
simplex $\sigma_x$ such that $x\in|\sigma_x|$, and moreover, $x$ lives
in the interior of $\sigma_x$. For each $g\in H$, there is exactly
one point $x_g \in \pi^{-1}(x)$ in $(g,\sigma_x)$.

Recall that, for every non-degenerate simplex $\sigma$ in $X$,
we have a continuous map $$g_{\sigma}: |\sigma| \to |X|$$ that is a
homeomorphism on the image when restricted to the interior. Moreover,
the images of these maps cover all $|X|$, and the topology of $|X|$ is
the one induced by this family of maps.

For any non-degenerate simplex $\sigma$ of dimension $m>n$, the set
$g_\sigma^{-1}(x)$ is nonempty if and only if $\sigma$ is incident to
$\sigma_x$. In that case, $g_\sigma^{-1}(x)$ is contained in the
boundary of $\sigma$. For each $k$-dimensional face of $\sigma$ (that is, a simplex
of the form $f=\partial_{i_1}\cdots\partial_{i_j}\sigma$), the
intersection of its interior with $g_{\sigma}^{-1}(x)$ is either the
interior of a polytope $P_f$ of dimension $m-j-n$ (if $f$ is a
degeneration of $\sigma_x$) or empty (in any other case). If we take a face of $f$, $f' = \partial_j f$ , the corresponding polytope $P_{f'}$ is exactly the intersection of $P$ with $f'$. That is, $g_\sigma^{-1}(x)$ is a polytopal complex $K_{x,\sigma}$, where the incidences between polytopes correspond to the incidences between the corresponding faces. Let $\{K_{x,\sigma}^1,\ldots,K_{x,\sigma}^{n_{x,\sigma}}\}$ be the connected components of $K_{x,\sigma}$. If we fix two points $y,y'\in K_{x,\sigma}^j$, belonging to faces $f,f'$ respectively, there must be a chain of faces$f=f_0,f_1,\ldots, f_s=f'$ such that each one is incident to the next. Now, if $y_h=x_g=y'_{h'}$ for some $h,h'\in H$, applying Lemma~\ref{lemma:caras-degeneradas-pegan} to the sequence, we obtain that $h=h'$. That is, the points on the same connected component that are mapped to $x_g$ have the same associated element of the group.

Finally, consider $U$ a neighborhood of $x$ that is small enough to
ensure that, for every non-degenerate simplex $\sigma$, $U \cap
|\sigma| = U_{x,\sigma}^1 \cup \cdots \cup U_{x,\sigma}^{n_{x,\sigma}}$,
where $U_{x,\sigma}^k$ is a neighborhood of
$K_{x,\sigma}^k$, and the unions are disjoint. By the previous
remark, there will be $h_{x,\sigma}^1,\ldots,h_{x_\sigma}^{n_{x,\sigma}} \in H$  such that the corresponding $(h_{x,\sigma}^k,U_{x,\sigma}^k)$ glued together form a neighborhood of $x_g$; and moreover, this gluing will mimic the one of $U_{x,\sigma}^k$ (because the simplicial map respects faces), so
this neighborhood will be homeomorphic to $U$.

\end{proof}

The previous results ensure us that we can construct covers of a simplicial set. To obtain the universal cover, we just need to choose $H$ and $\chi$ like this:

\begin{proposition}
If $\chi: \pi_1(|X|,v_*)\to H$ is an isomorphism, $|\tilde{X}|$ is the universal cover of $|X|$ and $\pi$ induces the corresponding covering map.
\end{proposition}
\begin{proof}
 We just have to show that the induced map in the fundamental group is trivial.
 An element of $\pi_1(|\tilde{X}|,(v_*,1))$ can be represented by a sequence of edges $(g_0,e_0),\ldots,(g_n,e_n)$. The formulas that define $\tilde X$ ensure that $g_{i+1}= g_i \cdot \tau (e_i) $. So, since the path must be closed, we get $\tau(e_0)\cdots \tau(e_n)=1$. Since this is exactly the element of $\pi_1(|X|,v_*)$ represented by the path $e_0,\ldots,e_n$, this means that the image of the path by $\pi$ is trivial in $\pi_1(|X|,v_*)$.
\end{proof}

In fact, it is easy to adapt the previous proof to see that, in general, we get the regular cover that corresponds to the group representation $\chi$.

These results lead to Algorithm~\ref{alg:universal_cover_simplicial_maximaltree}.

\begin{algorithm}[ht]
\label{alg:universal_cover_simplicial_maximaltree}
\caption{Cover of a surjective group morphism on the fundamental group}
\setcounter{AlgoLine}{0} 
\KwIn{
\begin{itemize}
    \item a connected simplicial set $X$,
    \item  a multiplicative group $H$,
    \item a surjective group morphism $\chi: \pi_1(|X|,v_*)\to H$.
\end{itemize}}
\KwOut{a simplicial set $\tilde{X}$ which is a cover of $X$.}
\lnl{p:maptau}
Define the map $\tau: X \rightarrow H$ by the recursive rules explained at the beginning of this section. \\
\lnl{p:simplices} Define the set of simplices of $\tilde{X}$ as $\tilde{X}_n = H \times X_n$.\\
\lnl{p:degeneracies} Define the degeneracy operators as $\tilde{s}_i(h,\sigma) = (h,s_i\sigma)$ for $h \in H $ and $\sigma \in X$. \\
\lnl{p:faces} Define the face maps, for an $n$-dimensional simplex $(h,\sigma)$, by $\tilde{\partial}_i(h,\sigma) = (h,\partial_i\sigma)$, for $i<n$, and $\tilde{\partial}_n(h,\sigma) = (h \cdot \tau(\sigma)^{-1}, \partial_n\sigma)$.\\
\Return $\tilde{X}$\\
\end{algorithm}

As mentioned above, if $H=\pi_1(|X|,v_*)$ (and $\chi$ is an isomorphism), then the output $\tilde{X}$ is a simplicial model for the universal cover of $X$.

If we want to work with simplicial sets of infinite type and implement our algorithm in Kenzo, then we can consider the effective homology method and the homological perturbation theorems as will be explained in the next section. To this aim, we observe that the map $\tau$ described before can also be seen as a twisting operator $\tau: X \rightarrow K(H,0)$, where $K(H,0)$ is the simplicial set defined as $K(H,0)_n=H$ for every $n$ and all face and degeneracy maps equals to the identity map (this simplicial set is an Eilenberg--MacLane space with $\pi_0(K(H,0),1)=H$ and $\pi_i(K(H,0),1)=0$ for all $i>0$). In fact, any map $\tau: X^{ND}_1 \rightarrow H$ (where  $X^{ND}_1$ is the set of non-degenerate edges of $X$) satisfying $\tau(\partial_1 x) = \tau(\partial_2x) \cdot \tau(\partial_0x) $ for all $x\in X^{ND}_2$ can be extended to a map $\tau : X \rightarrow K(H,0)$  as described at the beginning of this section (satisfying the properties of the twisting operator).
Then, the simplicial set $\tilde{X}$ is isomorphic to the twisted cartesian product defined by $\tau: X \rightarrow K(H,0)$, following Definition~\ref{def:twcrpr}. We obtain in this way Algorithm~\ref{alg:universal_cover_simplicial_product}.

\begin{algorithm}[ht]
\label{alg:universal_cover_simplicial_product}
\caption{Universal cover of a simplicial set as a twisted cartesian product}
\setcounter{AlgoLine}{0} 
\KwIn{
\begin{itemize}
\item a connected simplicial set $X$,
\item the fundamental group of $X$, $\pi_1(|X|,v_*) =H$,
\item a map defined on the non-degenerate edges of $X$, $\tau: X^{ND}_1 \rightarrow H$, such that $\tau(\partial_1 x) = \tau(\partial_2x) \cdot \tau(\partial_0x) $ for all $x\in X^{ND}_2$.
\end{itemize}}
\KwOut{a twisting operator $\tau : X \rightarrow K(H,0)$ that defines the simplicial model of the universal cover of $X$, expressed as a twisted cartesian product $K(H,0) \times_\tau X$.}
\lnl{p:maptauprime}
Define the map $\tau': X \rightarrow K(H,0)$ induced by  $\tau: X \rightarrow H$ (defined by the recursive rules explained at the beginning of this section) and the relation $K(H,0)_n = H$ for all $n \geq 0$. \\
\lnl{p:twistedproduct} Construct the simplicial set $K$ as the twisted cartesian product associated to $\tau'$ as described in Definition~\ref{def:twcrpr}.\\
\Return $K$\\
\end{algorithm}


\section{An algorithm computing the effective homology of the simplicial universal cover}
\label{sec:effective_homology_algorithm}

Using the simplicial version of the universal cover of a simplicial set as a twisted cartesian product  presented in Algorithm~\ref{alg:universal_cover_simplicial_product}, we can  determine the effective homology of universal covers by applying the homological perturbation theory as follows.

Let $X$ be a connected simplicial set (not necessarily of finite type) with fundamental group $\pi_1(|X|,v_*)=H$, and $\tau: X_1^{ND} \rightarrow H$ satisfying the hypotheses of Algorithm~\ref{alg:universal_cover_simplicial_product}. Let us also suppose now that $X$ has effective homology given by two reductions $C_*(X) \lrdc DX_*\rrdc EX_*$, where $EX_*$ is a chain complex of finite type (effective). We can consider the following diagram of reductions to determine the effective homology of the cartesian product of two simplicial sets~\citep{RS06}:

\begin{equation}
\xymatrix @R=6.0mm @C=2mm{
C_*(K(H,0)\times X) \ar @{=>>}^{\rho_1} [d] &  C_*(K(H,0))\otimes DX_* \ar @{=>>}_{\rho_2} [ld] \ar @{=>>}^{\rho_3} [rd]  \\
C_*(K(H,0))\otimes C_*(X) && C_*(K(H,0))\otimes EX_*
}
\label{eq:initial_reductions}
\end{equation}

The first reduction on the left, $\rho_1$, is the Eilenberg--Zilber reduction~\citep{EZ53}, defined from the chain complex associated with a cartesian product of two simplicial sets (in our case, $K(H,0)$ and $X$) to the tensor product of the chain complexes of the two factors. The reductions $\rho_2$ and $\rho_3$ are constructed, respectively, as the tensor product of the trivial reduction $C_*(K(H,0)) \rrdc C_*(K(H,0))$ with the two reductions of the effective homology of $X$, $DX_* \rrdc C_*(X)$ and $ DX_* \rrdc EX_*$ (following the formulas presented in Proposition~\ref{prop:red_tensor}).

Then, we try to apply the homological perturbation theory to these reductions, by considering the twisted cartesian product $K(H,0)\times_\tau X$, where $\tau$ is the twisting operator induced by the map $\tau: X_1^{ND} \rightarrow H$ as explained before. The differential of $C_*(K(H,0)\times_\tau X)$ can be seen as a perturbed version of the differential of $C_*(K(H,0)\times X)$, where the perturbation is given, for any $(k,x) \in K(H,0)_n \times X_n$, by
\begin{equation} \label{eq:pert_twist} \delta (k,x) \= (-1)^n \left[ (\tau (x) \cdot \partial_n k, \partial_n x) - (\partial_n k , \partial_n x) \right] .
\end{equation}

A major classical result, known as the twisted Eilenberg--Zilber theorem~\citep{Brown1959}, is obtained applying the Basic Perturbation Lemma to the Eilenberg--Zilber reduction, producing a new reduction $\hat{\rho}_1  : C_*(K(H,0)\times_\tau X) \rrdc C_*(K(H,0))\otimes_t C_*(X)$ where $C_*(K(H,0))\otimes_t C_*(X)$ is the perturbed chain complex obtained from $C_*(K(H,0))\otimes C_*(X)$ by introducing the perturbation induced (via the application of the BPL) by $\delta$. We denote this new perturbation by $\delta'$.
Now, using  the Trivial Perturbation Lemma (with the perturbation $\delta'$) we construct a reduction $\hat{\rho}_2 : C_*(K(H,0))\otimes_t DX_* \rrdc C_*(K(H,0))\otimes_t C_*(X)$ where $C_*(K(H,0))\otimes_t DX_*$ is a chain complex with the same underlying graded module as $C_*(K(H,0))\otimes DX_*$ and differential map obtained by a new perturbation $\delta''$. A last step is necessary for the construction of the effective homology of $K(H,0)\times_\tau X$: we need to apply the BPL to the reduction $\rho_3=(f_3,g_3,h_3): C_*(K(H,0))\otimes DX_* \rrdc C_*(K(H,0))\otimes EX_*$ (with the perturbation $\delta''$). This is the most difficult part of the process, since to apply the BPL one must verify that the nilpotency condition of the composition $h_3  \delta''$ is satisfied. In \citep{RS06}, a proof of this condition is given when the base space of the fibration is $1$-reduced, but this is not our case (the base space is $X$ and, in general, it is not simply connected; in fact, if $X$ is simply connected then it is its own universal cover). Therefore, a different proof of the nilpotency condition must be developed. To do this, we use the following results.

\begin{lemma}
\label{lemma:bottom_perturbation}
 Let $X$ be a simplicial set, $\pi_1(|X|,v_*)=H$ its fundamental group, $\tau: X^{ND}_1 \rightarrow H$ such that $\tau(\partial_1 x) = \tau(\partial_2x) \cdot \tau(\partial_0x) $  for all $x\in X^{ND}_2$, and $K(H,0)\times_\tau X$ as defined in Algorithm~\ref{alg:universal_cover_simplicial_product}. The twisted Eilenberg--Zilber reduction   $\hat{\rho}_1 = (\hat{f}_1,\hat{g}_1,\hat{h}_1) : C_*(K(H,0)\times_\tau X) \rrdc C_*(K(H,0))\otimes_t C_*(X)$ is in this case an isomorphism and the new perturbation $\delta'$ on $C_*(K(H,0))\otimes_t C_*(X)$ is given, for $(k \otimes x) \in (C_*(K(H,0))\otimes_t C_*(X))_n$, by:
 $$\delta'(k\otimes x) \= (-1)^n  [(\tau (x) \cdot  k) \otimes \partial_n x -  k \otimes \partial_n x].$$
\end{lemma}

\begin{proof}
Since $K(H,0)$ is defined as $K(H,0)_n = H$ for all $n \geq 0$ and all face and degeneracy maps are  identity maps, all simplices in $K(H,0)$ are degenerate except those in dimension $0$. Then, non-degenerate $n$-simplices of $K(H,0)\times_\tau X$ are pairs $(s_{n-1}\cdots s_0 k , x)$ with $k \in H$ and $x$ a non-degenerate $n$-simplex of $X$, and the generators of $C_*(K(H,0))\otimes_t C_*(X)$ are elements $k \otimes x$ with $k \in H$ and $x$ a non-degenerate $n$-simplex of $X$. In this way, the sets $(C_*(K(H,0))\otimes_t C_*(X))_n = C_0(K(H,0))\otimes C_n(X)$ and $C_n(K(H,0) \times X)$ are canonically isomorphic. It is also easy to observe that the differential map of both chain complexes  $C_*(K(H,0)\times_\tau X)$ and    $C_*(K(H,0))\otimes_t C_*(X)$ are the same, the maps $\hat{f}_1$ and $\hat{g}_1$ of the Eilenberg--Zilber reduction are in this case the canonical isomorphisms $\hat{f}_1(s_{n-1}\cdots s_0 k,x) = k\otimes x$ and $\hat{g}_1(k\otimes x) = (s_{n-1}\cdots s_0 k,x)$ for $k \in H$ and $x \in X_n$, and $\hat{h}_1$ is the null map. Finally, the perturbation $\delta'$ defined on  $C_*(K(H,0))\otimes C_*(X)$ is obtained by applying the BLP to the Eilenberg--Zilber reduction and the perturbation $\delta$, and is given by:

    \begin{equation*}
\begin{split}
\delta'(k\otimes x)  & = \hat{f}_1 \delta \hat{g}_1(k\otimes x) = \hat{f}_1 \delta (s_{n-1}\cdots s_0 k,x) \\
&  = \hat{f}_1((-1)^n  ((\tau (x) \cdot s_{n-2}\cdots s_0 k, \partial_n x) - ( s_{n-2}\cdots s_0 k , \partial_n x))) \\
& = (-1)^n  ((\tau (x) \cdot  k) \otimes \partial_n x -  k \otimes \partial_n x).
\end{split}
\end{equation*}

\end{proof}

\begin{lemma}
\label{lemma:top_perturbation}
Let $X$ be a simplicial set, $\pi_1(|X|,v_*)=H$ its fundamental group, $\tau: X^{ND}_1 \rightarrow H$ such that $\tau(\partial_1 x) = \tau(\partial_2x) \cdot \tau(\partial_0x) $ for all $x\in X^{ND}_2$, and $K(H,0)\times_\tau X$ as defined in Algorithm~\ref{alg:universal_cover_simplicial_product}. Let us suppose that $X$ has effective homology given by two reductions  $\rho^X_1=(f^X_1,g^X_1,h^X_1):DX_* \rrdc C_*(X)$ and $\rho^X_2=(f^X_2,g^X_2,h^X_2): DX_* \rrdc EX_*$. Then, the perturbation $\delta''$ induced on the chain complex $C_*(K(H,0))\otimes DX_*$ is given, for $(k \otimes x) \in (C_*(K(H,0))\otimes_t DX_*)_n$,  by:
$$\delta''(k\otimes x) = (-1)^n  [(\tau (f^X_1 (x)) \cdot  k) \otimes g^X_1 \partial_n f^X_1(x) -  k \otimes g^X_1 \partial_n f^X_1 (x)]$$
where $\partial_n:  C_n(X) \rightarrow C_{n-1}(X)$ and $\tau: C_n(X) \rightarrow H$ are obtained by extending by linearity the face $\partial_n$ of the simplicial set $X$ and the map $\tau: X \rightarrow H$.
\end{lemma}

\begin{proof}
It is directly obtained by applying the TPL to the reduction
$\rho_2: C_*(K(H,0))\otimes DX_* \rrdc C_*(K(H,0))\otimes C_*(X)$ and the formula for the perturbation $\delta'$ obtained in Lemma~\ref{lemma:bottom_perturbation}.
\end{proof}

\begin{lemma}
\label{lemma:composition}
Let $X$ be a simplicial set, $\pi_1(|X|,v_*)=H$ its fundamental group, $\tau: X^{ND}_1 \rightarrow H$ such that $\tau(\partial_1 x) = \tau(\partial_2x) \cdot \tau(\partial_0x) $ for all $x\in X^{ND}_2$, and $K(H,0)\times_\tau X$ as defined in Algorithm~\ref{alg:universal_cover_simplicial_product}. Let us suppose that $X$ has effective homology given by two reductions  $\rho^X_1=(f^X_1,g^X_1,h^X_1):DX_* \rrdc C_*(X)$ and $\rho^X_2=(f^X_2,g^X_2,h^X_2): DX_* \rrdc EX_*$. Let $\rho_3=(f_3,g_3,h_3)$ be the right reduction in diagram~(\ref{eq:initial_reductions}). Then, the composition  $h_3 \delta''$ is given, for $(k \otimes x) \in (C_*(K(H,0))\otimes_t DX_*)_n$, by:

$$h_3 \delta''(k\otimes x) = (-1)^n  [(\tau (f^X_1 (x)) \cdot  k) \otimes h^X_2 g^X_1 \partial_n f^X_1( x) -  k \otimes h^X_2 g^X_1 \partial_n f^X_1( x)].$$
\end{lemma}

\begin{proof}
Taking into account Proposition~\ref{prop:red_tensor} and the fact that in the trivial reduction $C_*(K(H,0)) \rrdc C_*(K(H,0))$ the map $h$ is null, we obtain that $h_3$ is defined,  for $(k \otimes x) \in (C_*(K(H,0))\otimes_t DX_*)_n$, as:
$$
h_3(k\otimes x) = (\id \otimes h^X_2)(k\otimes x) = k \otimes h^X_2(x)
$$
The result is easily obtained combining this formula with the expression of the perturbation $\delta''$ of Lemma~\ref{lemma:top_perturbation}.
\end{proof}

Let us recall that, to apply the BPL to the reduction $\rho_3$ in diagram~(\ref{eq:initial_reductions}), we need the composition $h_3 \delta''$ to be locally nilpotent, that is, for every $y \in C_*(K(H,0))\otimes_t DX_*$ there exists a nonnegative integer $m=m(y) \in \Nset$ such that
$(h_3 \delta'')^m(y)=0$. Once we have described the composition $h_3 \delta''$ in Lemma~\ref{lemma:composition}, we obtain the following result.

\begin{corollary}
\label{cor:composition_nilpotency}
Let $X$ be a simplicial set, $\pi_1(|X|,v_*)=H$ its fundamental group, $\tau: X^{ND}_1 \rightarrow H$ such that $\tau(\partial_1 x) = \tau(\partial_2x) \cdot \tau(\partial_0x) $ for all $x\in X^{ND}_2$, and $K(H,0)\times_\tau X$ as defined in Algorithm~\ref{alg:universal_cover_simplicial_product}. Let us suppose that $X$ has effective homology given by two reductions  $\rho^X_1=(f^X_1,g^X_1,h^X_1):DX_* \rrdc C_*(X)$ and $\rho^X_2=(f^X_2,g^X_2,h^X_2): DX_* \rrdc EX_*$ and assume that the composition $h^X_2 g^X_1 \partial_n f^X_1 $ satisfies that, for every element $y\in DX_n$, there exists some natural number $m$ with $(h^X_2 g^X_1 \partial_n f^X_1)^m(y) = 0 $. Then,  we can apply the BPL to the  reduction $\rho_3$ of diagram~(\ref{eq:initial_reductions}).
\end{corollary}

Combining all the previous results, we obtain the following diagram of reductions.

\begin{equation}
\label{eq:twisted_reductions}
\xymatrix @R=6.0mm @C=1mm{
C_*(K(H,0)\times_\tau X) \ar @{=>>}^{\hat{\rho}_1} [d] &  C_*(K(H,0))\otimes_t DX_* \ar @{=>>}_{\hat{\rho}_2} [ld] \ar @{=>>}^{\hat{\rho}_3} [rd]  \\
C_*(K(H,0))\otimes_t C_*(X) && C_*(K(H,0))\otimes_t EX_*
}
\end{equation}
where the reductions $\hat{\rho}_1$, $\hat{\rho}_2$, and $\hat{\rho}_3$ have been obtained by applying, respectively, the BPL, TPL and BPL to the reductions ${\rho}_1$, ${\rho}_2$, and ${\rho}_3$ of diagram~(\ref{eq:initial_reductions}).

Finally, we observe that, if the group $\pi_1(|X|,v_*)=H$ is finite, then the simplicial  group $K(H,0)$ is of finite type (effective) and therefore the chain complex $C_*(K(H,0))\otimes_t EX_*$ is also effective. In this way, if the conditions of Corollary~\ref{cor:composition_nilpotency} are satisfied, then the effective homology of the simplicial set $K(H,0)\times_{\tau} X$ (our simplicial model of the universal cover of $X$) can be obtained as the composition of the reductions of the previous diagram. We obtain in this way Algorithm~\ref{alg:effective_homology_univ_cover}.

\begin{algorithm} [ht]
\label{alg:effective_homology_univ_cover}
\caption{Effective homology of the universal cover}
\setcounter{AlgoLine}{0} 
\KwIn{
\begin{itemize}
\item a connected simplicial set $X$,
\item the fundamental group of $X$, $\pi_1(|X|,v_*)=H$, with a finite number of elements,
\item a map defined on the non-degenerate edges of $X$, $\tau: X^{ND}_1 \rightarrow H$, such that $\tau(\partial_1 x) = \tau(\partial_2x) \cdot \tau(\partial_0x) $ for all $x\in X^{ND}_2$,
\item the effective homology of $X$ given by two reductions $\rho^X_1=(f^X_1,g^X_1,h^X_1):DX_* \rrdc C_*(X)$ and $\rho^X_2=(f^X_2,g^X_2,h^X_2): DX_* \rrdc EX_*$ such that the composition $h^X_2 g^X_1 \partial_n f^X_1$ satisfies the nilpotency condition.
\end{itemize}}
\KwOut{the effective homology of the simplicial model of the universal cover of $X$, $K(H,0) \times_\tau X$.}
\lnl{p:twistedproductK} Construct the cartesian product $K(H,0) \times X$.\\
\lnl{p:EZreduction} Construct the Eilenberg--Zilber reduction $\rho_1: C_*(K(H,0) \times X) \rrdc C_*(K(H,0))  \otimes C_*(X)$, which in this case has been proven to be an isomorphism.\\
\lnl{p:tnprreductions} Construct the reductions $\rho_2 : C_*(K(H,0))  \otimes DX_* \rrdc C_*(K(H,0))  \otimes C_*(X) $ and $\rho_3 : C_*(K(H,0))  \otimes DX_* \rrdc C_*(K(H,0))  \otimes EX_*$ as the tensor product of the trivial reduction $C_*(K(H,0)) \rrdc C_*(K(H,0))$ with the two reductions of the effective homology of $X$, $DX_* \rrdc C_*(X)$ and $ DX_* \rrdc EX_*$ respectively (following the formulas presented in Proposition~\ref{prop:red_tensor}).\\
\lnl{p:rho1BPL} Apply the Basic Perturbation Lemma (Theorem~\ref{thm:bpl}) to the reduction $\rho_1$ with the perturbation $\delta$ corresponding to the twisting operator $\tau$, obtaining a new reduction $\rho'_1: C_*(K(H,0) \times_\tau X) \rrdc C_*(K(H,0))  \otimes_t C_*(X)  $ and a perturbation $\delta'$ on $ C_*(K(H,0))  \otimes C_*(X)$.\\
\lnl{p:rho2TPL} Apply the Trivial Perturbation Lemma (Theorem~\ref{thm:tpl}) to the reduction $\rho_2$ with the perturbation $\delta'$, obtaining a new reduction $\rho'_2 : C_*(K(H,0))  \otimes_t DX_* \rrdc C_*(K(H,0))  \otimes_t C_*(X) $ and a perturbation $\delta''$ on $ C_*(K(H,0))  \otimes_t DX_*$.\\
\lnl{p:rho3BPL} Apply the Basic Perturbation Lemma (Theorem~\ref{thm:bpl}) to the reduction $\rho_3$ with the perturbation $\delta''$, obtaining a new reduction $\rho'_3 : C_*(K(H,0))  \otimes_t DX_* \rrdc C_*(K(H,0))  \otimes_t EX_*$.\\
\Return Composition of the reductions $\rho'_1$, $\rho'_2$, and $\rho'_3$ (see diagram~(\ref{eq:twisted_reductions}))\\
\end{algorithm}

This algorithm can be combined with the Whitehead tower method for computing homotopy groups (implemented in Kenzo for $1$-reduced simplicial sets with effective homology and then enhanced in \citep{mathematics-2021} to deal with simply connected simplicial sets), taking into account that the universal cover $\tilde{X}$ satisfies $\pi_i(\tilde X) \cong \pi_i(X)$ for $i\geq 2$. Given a non simply connected simplicial set $X$, we determine its universal cover (which is simply connected) and its effective homology thanks to Algorithms~\ref{alg:universal_cover_simplicial_product} and \ref{alg:effective_homology_univ_cover}, and we apply the algorithm for computing its homotopy groups implemented in \citep{mathematics-2021}.

\begin{example}
\label{exm:S1WS2}
Let us see an example to showcase that the finiteness condition of the fundamental group is not superfluous.
Consider the simplicial set $X$ with the following non-degenerate simplices:
\begin{itemize}
\item A $0$-dimensional simplex $v$.
\item A $1$-dimensional simplex $e$, with $\partial_0(e)=\partial_1(e)=v$.
\item A $2$-dimensional simplex $t$, with $\partial_0(t)=\partial_1(t)=\partial_2(t)=s_0(v)$.
\end{itemize}

The corresponding topological space $|X|$ is $\mathbb{S}^1\vee \mathbb{S}^2$; therefore, the fundamental group is $\mathbb{Z}$, and the non-trivial homology groups are $H_1(X)\cong H_2(X)\cong \mathbb{Z}$. However, the universal cover $\tilde X$ is
a real line with one copy of $\mathbb{S}^2$ attached to each integer point. This space is simply connected, but $H_2(\tilde X)$ is the direct sum of infinitely many copies of $\mathbb{Z}$. Hence, $\tilde X$ cannot have effective homology, because it has a homology group that is not finitely generated.

\end{example}

\section{Implementation and didactic examples}
\label{sec:implementation_and_examples}

Algorithm~\ref{alg:universal_cover_simplicial_maximaltree} has been implemented in SageMath~\citep{sagemath}. It is valid only for simplicial sets of finite type (with a finite number of non-degenerate simplices) and finite fundamental group.
Our implementation is already available in \sagemath 10.0.

We illustrate its usage with an example. We start by creating a simplicial
set with finite fundamental group. In this case, we take the complex corresponding to the usual
presentation of the symmetric group on 3 elements, and take its product with the projective space (note that this simplicial set depends on the specific presentation of the initial group, not on the group itself).

\lstdefinelanguage{Sage}[]{Python}
{morekeywords={True,False,sage,singular},
deletekeywords=[2]{set},
literate={sage.}{sage.}6
{complex}{complex}8,
sensitive=true,
}

\lstset{frame=none,
          showtabs=False,
          showspaces=False,
          showstringspaces=False,
          commentstyle={\ttfamily\color{dredcolor}},
          keywordstyle={\ttfamily\color{dbluecolor}\bfseries},
          stringstyle ={\ttfamily\color{dgraycolor}\bfseries},
          language = Sage,
	  basicstyle={\small \ttfamily},
	  aboveskip=.3em,
	  belowskip=.1em
          }

\begin{lstlisting}
sage: G = SymmetricGroup(3).as_finitely_presented_group()
sage: G
Finitely presented group < a, b | b^2, a^3, (a*b)^2 >
sage: C = simplicial_sets.PresentationComplex(G)
sage: RP3 = simplicial_sets.RealProjectiveSpace(3)
sage: S = C.product(RP3) ; S
Simplicial set with 12 non-degenerate simplices x
RP^3
sage: S.fundamental_group()
Finitely presented group < e0, e5, e9 | e0^2, e9^2,
e5^3, e0*e5*e0*e5^-1, (e5*e9)^2, (e9*e0)^2,
(e5*e0*e9)^2 >
sage: S.fundamental_group().cardinality()
12
\end{lstlisting}

Creating its universal cover takes 22 seconds in an Intel Core i7-10700, using 260MB of RAM:

\begin{lstlisting}
sage: SC = S.universal_cover()
sage: SC
Simplicial set with 4176 non-degenerate simplices
\end{lstlisting}

We can check that it is indeed simply connected:

\begin{lstlisting}
sage: SC.fundamental_group()
Finitely presented group <  |  >
\end{lstlisting}

And now we can compute its usual topological invariants, as any other simplicial set,
and compare them with the base space (this computation takes about a minute to complete):

\begin{lstlisting}
sage: [SC.homology(i) for i in range(6)]
[0, 0, Z^11, Z, 0, Z^11]
sage: [S.homology(i) for i in range(6)]
[0, C2 x C2, Z x C2, Z x C2 x C2, C2, Z]
\end{lstlisting}

Using the interface with \kenzo, we can apply the Whitehead tower method to compute also the higher homotopy groups of the cover (which coincide with the corresponding homotopy groups of the base space):

\begin{lstlisting}
sage: from sage.interfaces.kenzo import KFiniteSimplicialSet
sage: KSC = KFiniteSimplicialSet(SC)
sage: KSC.homotopy_group(2)
Multiplicative Abelian group isomorphic to Z x Z x Z x Z x Z x Z x
Z x Z x Z x Z x Z
\end{lstlisting}

Algorithms~\ref{alg:universal_cover_simplicial_product} and~\ref{alg:effective_homology_univ_cover} have been implemented as functions in the Kenzo system (the code can be found at~\citep{Kenzo-implementation}). As said before, one of the main advantages of the use of this system is that it allows one to work with spaces of infinite type and, moreover, it allows one to use the effective homology theory. In this way, it is possible to determine a simplicial model of the universal cover of simplicial sets of infinite type (with finite fundamental group) and determine its homology and homotopy groups.

To illustrate our programs and the power of the effective homology theory in our problem, we consider as a didactic example the following simplicial set of infinite type: we build in Kenzo the cartesian product of the projective plane with the ``semiline'' divided in intervals. The projective plane $\Rset P $ is given by non-degenerate simplices $\Rset P_0=\{v\}$, $\Rset P^{ND}_1=\{a\}$ and $\Rset P^{ND}_2=\{t\}$, and faces $\partial_0 a = \partial_1 a = v$, $\partial_0 t = \partial_2 t = a $ and $\partial_1 t = s_0 v$; it is a simplicial set of finite type. The ``semiline'' is represented as a simplicial set $A$ with non-degenerate simplices given by $A_ 0 = \{ n | n \in \Nset\}$, $A^{ND}_ 1 = \{ [n,n+1] | n \in \Nset\}$ and $\partial_0([n,n+1]) = n+1$, $\partial_1([n,n+1]) = n$. The simplicial set $A$ has an infinite number of non-degenerate simplices, but we can construct its effective homology in an explicit way as follows. We construct a reduction $\rho^A_2$:

$$
 \xymatrix {
 C_\ast(A) \ar @(l,u) []^{h^A_2} \ar @/^/[rr]^{f^A_2} & & C_*(*)\ar @/^/ [ll]^{g^A_2}
}
$$
where $C_*(*)$ is a chain complex with only one generator $*$ in degree $0$, $f^A_2$ is defined by $f^A_2(x)=*$ if $x\in A_0$ and $f^A_2(x)=0$ for all $x \in C_n(A)$ with $n > 0$, $g^A_2$ is given by $g^A_2(*)=0 \in A_0$ and $h^A_2(n)= [0,1]+[1,2]+\cdots + [n-1,n]$ (and $h^A_2(x) = 0$ for $x \in C_n(A)$ with $n > 0$). It is easy to verify that these maps satisfy the properties of reduction. The left reduction in the effective homology of $A_*$ is the identity reduction $ \rho^A_1 = \id: C_*(A) \lrdc C_*(X) $.

Now, the cartesian product $X= \Rset P \times A$ is also a simplicial set with effective homology (it is built automatically by Kenzo with the same idea as the reductions in diagram~(\ref{eq:initial_reductions})). Moreover, since $A$ is contractible, its fundamental group is equal to the fundamental group of $\Rset P$, that is $\pi_1(\Rset P)= \Zset / 2 \Zset$.

To apply Algorithm~\ref{alg:effective_homology_univ_cover} and construct the universal cover of $X$ with effective homology, it is also necessary to define  a map $\tau: X^{ND}_1 \rightarrow \Zset /2 \Zset$, and we can do it in the following way: let $e$ be an edge in $X^{ND}_1$, $e = (y,z)$ with $y \in \Rset P$ and $z \in A$. If $y = a$ (the unique edge in $\Rset P$), then we define $\tau(e)=1$; if $y \neq a$, then we define $\tau(e)=0$. Finally, it is not difficult to verify that the effective homology of $X$ satisfies the condition of Algorithm~\ref{alg:effective_homology_univ_cover}, that is, the composition $h^X_2 g^X_1 \partial_n f^X_1$ satisfies the nilpotency condition.

Applying now Algorithm~\ref{alg:effective_homology_univ_cover}, we can construct in Kenzo the simplicial model for the universal cover of $X$ with its effective homology as follows. To this aim, we build the cartesian product of the projective plane and the semiline (we omit the construction of the simplicial sets \texttt{proj-plane} and \texttt{semiline}) and we store it in the variable \texttt{X}. Then, we define the map $\tau: X^{ND}_1 \rightarrow \Zset /2 \Zset$ as a function that receives an edge and returns an element of the group (in this case, $0$ or $1$) and we store it in \texttt{X-twop-edges}. Finally, we call the function \texttt{universal-cover}:

\lstset{
  language=Lisp,
  basicstyle={\small \ttfamily},
}

\begin{lstlisting}
> (setf X (crts-prdc proj-plane semiline))
[K21 Simplicial-Set]
> (setf X-twop-edges
  #'(lambda (edge)
      (with-crpr (dgop1 gmsm1 dgop2 gmsm2) edge
        (if (and (eq 0 dgop1) (eq gmsm1 'a))
            1 0))))
#<Interpreted Function (unnamed) @ #x20e50d82>
> (setf X-univ-cover (universal-cover X (cyclicgroup 2) X-twop-edges))
[K45 Simplicial-Set]
\end{lstlisting}

As said before, this allows us to determine its homology and homotopy groups. For instance, we compute the homotopy group of dimension $5$:

\begin{lstlisting}
> (homotopy X-univ-cover 5)
Homotopy in dimension 5 :
Component Z/2Z
\end{lstlisting}
which indeed is the correct result, since $\pi_5(\mathbb{S}^2)=\mathbb{Z}/2\mathbb{Z}$.

\section{Study of some algebraic topology constructors and tools}
\label{sec:constructors_and_tools}

\subsection{Cartesian products}

Given two connected simplicial sets $X$ and $Y$, the universal cover of the cartesian product $X\times Y$ is well known to be isomorphic to the cartesian product of the universal covers of the factors. If $X$ and $Y$ satisfy the hypotheses of Algorithm~\ref{alg:effective_homology_univ_cover}, then we can compute the effective homology of the universal covers $\Tilde{X}$ and $\Tilde{Y}$. Applying then the algorithm for computing the effective homology of a cartesian product, we would obtain the effective homology of the universal cover of $X\times Y$.

\subsection{Discrete vector fields}

Discrete vector fields (an important idea of Robin Forman's Discrete Morse Theory~\citep{For98}) are a tool allowing one to obtain the effective homology of complicated spaces in some interesting situations. For example, a discrete vector field can be defined producing the well-known  Eilenberg--Zilber theorem~\citep{EZ53}, which, as seen before, provides a reduction from the chain complex of the cartesian product of two simplicial sets to the tensor product of the chain complexes of the factors. Moreover, the same vector field is valid for the  case of twisted cartesian products, described in the twisted Eilenberg--Zilber theorem~\citep{Brown1959}. See~\citep{RS10} for the definition of these vector fields and for details of the following definitions.

Let $C_\ast=(C_n,d_n)$ be a free chain complex with distinguished $\Z$-bases $B_n \subseteq C_n$, whose elements we call $n$-\emph{cells}.

\begin{definition}
\label{defn:vector_field}
A \emph {discrete vector field} $V$ on $C_\ast$ is a collection of pairs of cells $V = \{(x_k; y_k)\}_{k\in K}$ such that:
\begin{itemize}
\item Every $x_k$ is an element of some $B_n$, in which case $y_k\in B_{n+1}$. The degree $n$ depends on $k$ and in general is not constant.
\item Each component $x_k \in B_n$ is a \emph{regular face} of the corresponding $y_k \in B_{n+1}$ (that is, the coefficient of $x_k$ in $d_{n+1}(y_k)$ is $+1$ or $-1$).
\item Each generator (cell) of $C_\ast$ appears at most once in $V$.
\end{itemize}
\end{definition}

\begin{definition}
A pair $(x_j ;y_j )$ of $V$ is called a \emph{vector}. 
The cells $x_j$ and $y_j$ are called respectively a \emph{source cell} and a \emph{target cell}.
A cell $x\in B_n$ that does not appear in the discrete vector field $V$ is called a \emph{critical cell}.
\end{definition}

\begin{definition}
Given a discrete vector field $V$, a \emph{$V$-path $\pi$ of
degree $n$ and length $m$} is a sequence $\pi =\{(x_{j_k}; y_{j_k})\}_{0 \leq k <
m}$ such that:
\begin{itemize}
\item
Every pair $(x_{j_k}; y_{j_k})$ is a vector of
$V$ and $y_{j_k}$ is an $n$-cell.
\item
For every $0 < k < m$, the component $x_{j_k}$ is a \emph{face} of
$y_{j_{k-1}}$ (meaning that the coefficient of $x_{j_k}$ in $d_{n}(y_{j_{k-1}})$ is non-null), non necessarily regular but different from
$x_{j_{k-1}}$.
\end{itemize}
\end{definition}

\begin{definition}
A discrete vector field $V$ is called \emph{admissible} if, for
every $n \in \Zset$, a function $\lambda_n: B_n \to \Nset$ is
provided such that the length of every $V$-path starting from
$x \in B_n$ is bounded by $\lambda_n(x)$.
\end{definition}

The following result, due to Forman \cite[\S~8]{For98}, has been generalized in~\citep{RS10} to the case of chain complexes not necessarily of finite type.

\begin{theorem}(\cite{For98,RS10})
\label{thm:dfv-reduction}
Let $C_\ast=(C_n,d_n,B_n)$ be
a free chain complex and $V = \{(x_k; y_k)\}_{k\in K}$ be an admissible discrete vector
field on $C_\ast$. Then the vector field $V$ defines a canonical reduction $\rho = (f, g, h) :
(C_n, d_n)\rrdc (C^c_n, d'_n)$ where $C^c_n$ is the free $\Zset$-module generated by the critical
$n$-cells and $d'_n$ is an appropriate differential canonically defined from $C_\ast$ and $V$.
\end{theorem}

Let us suppose now that $X$ is a simplicial set, $\pi_1(|X|,v_*)=H$ its fundamental group and $V$ a vector field defined on the associated normalized chain complex, $C^N_*(X)$ (the $n$-cells are in this case the non-degenerate $n$-simplices of $X$). The vector field $V$ induces a vector field $\Tilde{V}$ on the (normalized chain complex of the) universal cover $\Tilde{X}$ taking as vectors of $\Tilde{V}$, for each element $h \in H)$:
\begin{itemize}
    \item
the pairs $((h,x),(h,y))$ for each vector $(x,y)$ of $V$ such that $x \in X_n$ and $x \neq \partial_{n+1}y$
\item and the pairs $((h,x),(h \cdot \tau(y),y))$ for each vector $(x,y)$ in $V$ such that $x \in X_n$ and $x = \partial_{n+1}y$.
\end{itemize}

\begin{proposition}
If $V$ is an admissible vector field on a simplicial set $X$, then $\Tilde{V}$ is an admissible vector field on the universal cover $\Tilde{X}$.
\end{proposition}

\begin{proof}
First of all, we prove that $\Tilde{V}$ satisfies the  properties of Definition~\ref{defn:vector_field}. The first and third one are trivial. For the second one, we have two cases:
\begin{itemize}
    \item Given a vector $(x,y)$ of $V$ such that $x \in X_n$ and $x \neq \partial_{n+1}y$, one has that $x = \partial_jy$ for $0 \leq j < n+1$. Then, the corresponding vectors in $\Tilde{V}$, that is, the pairs $((h,x),(h,y))$ for $h \in H$  satisfy $\partial_j(h,y) = (h,x)$ and the coefficient of $(h,x)$ in $\Tilde{d}_{n+1}(h,y)$ is equal to the coefficient of $x$ in $d_{n+1}(y)$, so that $(h,x)$ is a regular face of its associated target cell $(h,y)$.
    \item Given a vector $(x,y)$ of $V$ such that $x \in X_n$ and $x = \partial_{n+1}y$, then, for all  $h \in H$, $\partial_{n+1}(h \cdot \tau(y),y) = (h \cdot \tau(y) \cdot \tau^{-1}(y),\partial_{n+1} y ) = (h,x)$. Moreover, if there are other faces $\partial_j(h \cdot \tau(y),y) = (h,x)$ with $0 \leq j < n+1$, then $\partial_j y  = x$, which implies that the coefficient of  $(h,x)$ in $\Tilde{d}_{n+1}(h \cdot \tau(y),y)$ is again $+1$ or $-1$ and therefore $(h,x)$ is a regular face of the associated target simplex $(h \cdot \tau(y),y)$.
\end{itemize}

To prove that $\Tilde{V}$ is admissible, it suffices to define $\Tilde{\lambda}_n: \Tilde{X}_n \rightarrow \Nset$ as $\Tilde{\lambda}_n(h,x) = \lambda_n(x)$, where $\lambda_n$ is the function proving the admissibility of the vector field $V$.

\end{proof}

\begin{corollary}
\label{cor:dvf_universal_cover}
    Let $X$ be a connected simplicial set, $\pi_1(|X|,v_*)=H$ its fundamental group, $\tau: X^{ND}_1 \rightarrow H$ such that $\tau(\partial_1 x) = \tau(\partial_2x) \cdot \tau(\partial_0x) $ for all $x\in X^{ND}_2$, and $K(H,0)\times_\tau X$ as defined in Algorithm~\ref{alg:universal_cover_simplicial_product}. Let $V$ be an admissible discrete vector field on $X$ such that the set $X^c_n$ of $n$-critical cells is finite for every dimension $n$. If $H$ is finite, a reduction  $\rho: C_*(\Tilde{X}) \rightarrow C_*(\Tilde{X}^c)$ can be built, where  $C_*(\Tilde{X}^c)$ is an effective chain complex.
\end{corollary}

A particular case of this situation is obtained for the Eilenberg--MacLane spaces $K(\Zset,1)$ and $K(\Zset/m \Zset,1)$, for $m \geq 2$. In both cases, admissible discrete vector fields can be defined on the associated chain complexes, producing a reduction to an effective chain complex (see~\citep{RS10} and~\citep{Rom10} respectively). Thanks to our previous results, these vector fields can be used to define vector fields on the associated universal covers.  For $K(\Zset/m \Zset,1)$, this provides the effective homology of its universal cover. The case of  $K(\Zset,1)$ (more concretely, the case of infinite fundamental group) will be studied in Section~\ref{sec:infinite_fundamental_group}.

Discrete vector fields can also be used to compute the effective homology of twisted cartesian products~\citep{RS10} and classifying spaces~\citep{DRRS22}. If the initial spaces are effective (of finite type), then these vector fields combined with our Corollary~\ref{cor:dvf_universal_cover} provide the effective homology of the universal cover. If some of the factor of the cartesian product or the initial space in the classifying space is of infinite type, then other reductions appear in the effective homology of the new spaces and these reductions should be studied to see if they satisfy the hypotheses of Algorithm~\ref{alg:effective_homology_univ_cover}.

\section{The case of abelian infinite fundamental group}
\label{sec:infinite_fundamental_group}

As we saw in Example~\ref{exm:S1WS2}, when the fundamental group is not finite,
the universal cover might not have effective homology. However, in some cases its homology groups can be computed nevertheless.

As before, let $X$ be a connected simplicial set, and $G$ its fundamental group. For any $n\in \mathbb{N}$, the chain group $C_n(X)$ is freely generated by the set $X^{ND}_n$ of non-degenerate simplices of dimension $n$. By construction of $\tilde X$, $C_n(\tilde X)$ is freely generated by $G\times X^{ND}_n$. That is, the abelian group structure of $C_n(\tilde X)$ is isomorphic to $\mathbb{Z}[G]\otimes C_n(X)$.

Moreover, it is easy to check that the boundary map

$$
\begin{array}{rccc}
\tilde d_n: & C_n(\tilde{X}) & \longrightarrow & C_{n-1}(\tilde{X}) \\
& (g,\sigma) & \mapsto & \sum_{i=0}^{n-1}(-1)^i(g,\partial_i\sigma) + (-1)^n (g\cdot \tau(\sigma)^{-1},\partial_n\sigma)
\end{array}
$$
is a morphism of $\mathbb{Z}[G]$-modules. So, the homology groups of $\tilde X$ will be, in fact, the homology groups of the chain complex of $\mathbb{Z}[G]$-modules

$$
\cdots \xleftarrow{\tilde{d}_{n-1}} C_{n-1}(\tilde X) \xleftarrow{\tilde{d}_n}  C_{n}(\tilde X) \xleftarrow{\tilde{d}_{n+1}}  C_{n+1}(\tilde X) \xleftarrow{\tilde{d}_{n+2}} \cdots
$$

If the simplicial set $X$ is effective itself, this chain complex can be computed explicitly; and if the group $G$ is abelian, the group algebra $\mathbb{Z}[G]$ is isomorphic to a quotient of a polynomial ring. In this case, a generating set of $\ker(\tilde{d}_n)$ will be given by the syzygies of the images of the generators (these syzygies can be computed using Gröbner basis methods over modules). The module $\ker(\tilde{d}_n)$ will be a quotient of the free module generated by those generators, where the relations are again the syzygies between the generators.
Adding the images of the generators of $C_{n+1}(\tilde X)$, expressed in the generators of $\ker(\tilde d_n)$, we obtain presentation matrices for $H_n(\tilde X)$ as $\mathbb{Z}[G]$-modules.

\begin{algorithm}[ht]
\label{alg:twisted_homology}
\caption{Twisted homology of an effective simplicial set.}
\setcounter{AlgoLine}{0} 
\KwIn{
\begin{itemize}
    \item an effective connected simplicial set $X$ with abelian fundamental group $G$,
    \item  a positive integer $n$,
\end{itemize}}
\KwOut{a presentation matrix of $H_n(\tilde X)$ as $\mathbb{Z}[G]$-module.}

\lnl{p:grpalg}
Express the group $G$ as direct sum of cyclic groups: $G = \mathbb{Z}/t_1\mathbb{Z}\oplus\cdots\oplus \mathbb{Z}/t_l\mathbb{Z}\oplus \mathbb{Z}^f$\\
\lnl{p:plnring}
Take the polynomial ring $R = \mathbb{Z}[x_1,\bar x_1,\ldots,x_{l+f},\bar x_{l+f}]$ and the ideal
$I = (x_1\bar x_1 -1, \ldots , x_{l+f}\bar x_{l+f}-1, x_{1}^{t_1}-1,\ldots, x_l^{t_l}-1)$.
Consider the isomorphism $\nu:\mathbb{Z}[G]\to R/I.$
\\
\lnl{p:lstbs}
Construct ordered lists $L_{n-1}=(\sigma^{n-1}_1,\ldots\sigma^{n-1}_{m_{n-1}}), L_n = (\sigma^{n}_1,\ldots\sigma^{n}_{m_{n}}),L_{n+1}=(\sigma^{n-1}_1,\ldots\sigma^{n+1}_{m_{n+1}})$ with the non-degenerate simplices of dimensions $n-1, n$ and $n+1$.\\
\lnl{p:mtrcs}
Construct a $m_{n-1}\times m_n$ matrix $M^n$, and a $m_n\times m_{n+1}$ matrix $M^{n+1}$, both initialized with zero entries.\\
\lnl{p:fillmatrcs} For each $\sigma^n_j$, and for each $0\leq d< n$, if
$\partial_d(\sigma^n_j)=\sigma^{n-1}_i$, increase $M^n_{i,j}$ by $(-1)^d$. If
$\partial_n(\sigma^n_j)=\sigma^{n-1}_k$, increase $M^n_{k,j}$ by
$(-1)^n\nu(\tau(\sigma^n_j)^{-1})$ .\\
 \lnl{p:fillmatrcs2} For each $\sigma^{n+1}_j$, and for each $0\leq d< n+1$, if $\partial_d(\sigma^{n+1}_j)=\sigma^{n}_i$, increase $M^{n+1}_{i,j}$ by $(-1)^d$. If $\partial_{n+1}(\sigma^{n+1}_j)=\sigma^{n}_k$, increase $M^{n+1}_{k,j}$ by $(-1)^n\nu(\tau(\sigma^n_j)^{-1})$ .\\
\lnl{p:syz1} Compute the syzygy matrix $N$ between the columns of $M^n$ (that is,
the columns of $N$ generate the right kernel of $M^n$.\\
\lnl{p:relsexp} For each column $c_i$ of $M^{n+1}$, express $c_i$ as a linear combination of the rows of $N$ obtaining a matrix $R_1$ such that $M^{n+1} = N R_1$.\\
\lnl{p:reslsint} Compute the syzygy matrix between the columns of $N$, obtaining a matrix $R_2$ whose columns generate the right kernel of $N$.\\
\lnl{p:outrels} Output the matrix obtained by stacking $R_1^T$ and $R_2^T$.
\end{algorithm}

We implemented this approach (see~\citep{sagemath-pr-twistedhomology}) and it will be included in future releases of \sagemath.

Here we present an example that showcases how this approach can handle Example~\ref{exm:S1WS2}:

\lstset{
  language=Sage,
  basicstyle={\small \ttfamily},
}

\begin{lstlisting}
sage: X = simplicial_sets.Sphere(1).wedge(simplicial_sets.Sphere(2))
sage: X.twisted_homology(1)
Quotient module by Submodule of Ambient free module of rank 0
over the integral domain Multivariate Polynomial Ring in f1, f1inv
over Integer Ring
Generated by the rows of the matrix:
[]
sage: X.twisted_homology(2)
Quotient module by Submodule of Ambient free module of rank 1
over the integral domain Multivariate Polynomial Ring in f1, f1inv
over Integer Ring
Generated by the rows of the matrix:
[f1*f1inv - 1]
\end{lstlisting}

Notice how $H_1(\tilde X)$ is trivial (it lives inside a module of rank $0$), whereas
$H_2(\tilde X)$ is isomorphic to $\mathbb{Z}[t,t^{-1}]$; that is, as an abelian group, it is freely generated by one generator for each integer.

\section{Conclusions and further work}
\label{sec:conclusions}

This paper presents contributions in order to perform new computations in algebraic topology. So far, SageMath allows one to compute homology groups of chain complexes (and simplicial sets) of finite type. Furthermore, Kenzo (and its interface in SageMath) allows one to compute homology and homotopy groups for certain types of infinite  spaces by means of effective homology, being the only software able to perform this kind of computation. The computation of homotopy groups requires the simplicial sets to be simply connected.   This work aims to generalize the computational capability of both programs, allowing one to perform such computations involving simplicial sets that are not simply connected.

For this task, we rely on the idea of universal cover of a space, which is another simply connected space, with a covering map. This universal cover satisfies that the higher homotopy groups coincide with those of the initial space. For this purpose, we have first developed a simplicial version of the universal cover, which is represented by means of twisted cartesian products. Then, we have developed an algorithm for computing the effective homology of the simplicial universal cover, based on the application of the homological perturbation theory. The algorithms have been implemented in both SageMath and Kenzo.

As a future work, our intention is to take advantage of the existing interface between SageMath and Kenzo~\citep{mathematics-2021}, so that all the calculations that Kenzo is able to perform thanks to this work (mainly, homotopy group calculations through the universal cover), can also be performed directly in SageMath in a transparent manner. For this case, the connection between both programs is not trivial and requires a deep study of the possibilities: to build the universal cover in Kenzo, we need a function from the edges to the fundamental group. This would require to \emph{somehow} input functions to Kenzo from SageMath code.

%
%


\end{document}